\documentclass[12pt]{amsart}
\usepackage{epsfig,amssymb,epic,eepic,color}
\usepackage[all]{xy}
 \usepackage[T1]{fontenc}        
\newtheorem{thm}{Theorem}[section]

\newtheorem{prop}[thm]{Proposition}

\newtheorem{lemme}[thm]{Lemma}
\newtheorem{cor}[thm]{Corollary}

\newtheorem{rien}[thm]{}
 
\newcommand{\be}{\begin{enumerate}}
\newcommand{\ee}{\end{enumerate}}
\newcommand{\bi}{\begin{itemize}}
\newcommand{\ei}{\end{itemize}}

\def\R{\mathbb{R}}

\def\Z{\mathbb{Z}}

\def\Om{\Omega}

\def\ga{\gamma}    
\def\Ga{\Gamma}

\def\al{\alpha}
\def\be{\beta}

\def\De{\Delta}
\def\vp{\varphi}

\def\la{\lambda}

\def\si{\sigma}
\def\Si{\Sigma}

\def\ep{\varepsilon}

\def\ds{\displaystyle}

\def\nd{\noindent}
\def\bull{\hfill$\Box$}
\def\proof{\nd {\bf Proof.\ }}

\begin{document}
\vskip 1cm
%\today
\begin{center}
{\sc Regularization of
$\Gamma_1$-structures in dimension 3}
\end{center}

\title{}
\author{ Fran\c cois Laudenbach}
\address{Universit\'e de Nantes, UMR 6629 du CNRS, 44322 Nantes, France}
\email{francois.laudenbach@univ-nantes.fr}

\author{ Ga\"el Meigniez}
\address{Universit\'e de Bretagne-Sud,  L.M.A.M., 
BP 573, F-56017 Vannes, France}
\email{Gael.Meigniez@univ-ubs.fr}

%\date{This version: April 2, 2009.}



\thanks{FL supported by ANR Floer Power}

\begin{abstract} For %co-orientable
 $\Gamma_1$-structures on 3-manifolds, we give a very simple proof of Thurston's regularization theorem, first proved in \cite{thurston}, without using  %the perfectness of the group of diffeomorphisms of the circle, nor 
 Mather's homology equivalence. % between BDiff$_c(\R)$ and $\Om$ B$(\Ga_1)_+$.
%Diff$^\infty(S^1)$.
% In particular it applies in  classes of regularity where perfectness does not hold (or is %not known to hold). 
Moreover, in the co-orientable case, 
 the resulting foliation can be chosen of a precise kind, namely
an ``open book foliation modified by suspension''. There is also a model in the non co-orientable case.

\end{abstract}

\maketitle

\thispagestyle{empty}
\vskip 1cm
\section{Introduction}
\medskip

A $\Gamma_1$-structure $\xi$, in the sense of A. Haefliger,
 on a manifold $M$ is given by
a line bundle $\nu=(E\to M)$, called the {\it normal bundle} to $\xi$, and
 a germ of codimension-one foliation $\mathcal F$
along the zero section, which is required to be transverse to the fibers
 (see \cite{haef}). 
To fix ideas, consider the  co-orientable case, that is, 
the normal bundle is trivial: $E\cong M\times\R$; for the general case see
 section 
\ref{comments}.
The $\Ga_1$-structure $\xi$ is said to 
be {\it regular} when the foliation $\mathcal F$ is transverse to the
 zero section, in which
case the pullback of $\mathcal F$ to $M$  
is a genuine foliation on $M$. 
A homotopy of $\xi$ is defined as a  
$\Gamma_1$-structure on $M\times [0,1]$ inducing $\xi$ on $M\times\{0\}$. A 
{\it regularization theorem} should claim that any $\Ga_1$-structure is 
homotopic to a regular one. It is not true in general. 
An obvious necessary condition is
that 
$\nu$ must embed into the tangent bundle $\tau M$. When $\nu$ is trivial and 
$\dim M=3$ this condition is fulfilled. 

%Usually, when nothing is specified, it is meant that 
The $C^\infty$ category is understood in the sequel, 
unless otherwise specified. In particular $M$ is $C^\infty$. 
One calls $\xi$ a $\Ga^r_1$-structure ($r\geq 1$) if it is 
tangentially $C^\infty$ and transversely $C^r$, 
that is, the foliation charts are 
$C^r$ in the direction transverse to the leaves.
We will prove the following theorem.

\begin{thm} \label{reg}If $M$ is a closed 3-manifold and $\xi $
 a $\Ga_1^r$-structure, 
$r\geq 1$,
 whose normal bundle is trivial, then $\xi$ is homotopic to a regular 
 $\Ga_1^r$-structure.
\end{thm}

Moreover, the resulting foliation of $M$ may have its tangent plane field 
in a prescribed
homotopy class (see proposition \ref{nu}).

 %When $r>2$ 
 This theorem is  a particular case of a general regularization 
theorem due to W. Thurston (see \cite{thurston}). Thurston's proof was based on the deep result due to J. Mather \cite{mather0}, \cite{mather1}:
the homology equivalence between the classifying space of the  group Diff$_c(\R)$ 
endowed with the discrete topology and the loop space $\Om$B$(\Ga_1)_+$.
%; and also, in dimension 3, 
%the famous perfectness of the {\it discrete} group Diff$^r(S^1)$, due to 
%M. Herman when $r=\infty$ and to J. Mather when $r>2$
 %(see \cite{herman, herman2}
%\cite{mather,mather2}), and which is unknown at the present time when $r=2$. 
We present  a proof of this regularization theorem which does not need this result.
%As \'Etienne Ghys pointed out to us, the proof also works  in the class of
% regularity
%$r=1+bv$ (see proposition \ref{bv}), the class he calls the Denjoy class 
%(as it is the class of regularity where 
%Denjoy's conjugation theorem holds). Here Diff$^{1+bv}(S^1)$ means the set of 
%$C^1$-diffeomorphisms  $f$ 
%whose derivative has a bounded variation, that is, $df$ (or $\log(df)$) 
%s the primitive of a Radon measure. J. Mather proved in \cite{mather3} 
%that the group Diff$^{1+bv}(S^1)$ is not perfect: there is a 
%surjective homomorphism Diff$^{1+bv}(S^1)\to\R$. 
A regularization theorem in all dimensions, still
avoiding any difficult result,
is provided in 
 \cite{meigniez}. But there are reasons 
 for considering  the dimension 3 separately. 
 
 Our proof
% not only solves a regularity problem but it 
provides models realizing each homotopy class of $\Ga_1$-structure.
%a model for each homotopy class of $\Ga_1$-structure. 
The models are based on the notion of
 {\it open book decomposition}. 
 Recall that such a structure on $M$ consists of a link $B$ in $M$,
 called the {\it binding},
 and a fibration $p: M\smallsetminus B\to S^1$ such that, 
for every $\theta\in S^1$, 
 $p^{-1}(\theta)$ is the interior of an embedded surface, 
called the {\it page} $P_\theta$, 
 whose boundary is the binding. The existence of open book decomposition 
could be 
 proved by J. Alexander when $M$ is orientable, as a consequence
 of \cite{alex} (every orientable closed 3-manifold is a branched cover of 
the 3-sphere) and \cite{alex2} 
 (every link can be braided); but he was ignoring this concept which was 
introduced by 
 H. Winkelnkemper in 1973 \cite{wink}. Henceforth,  we refer to the more 
flexible construction by E. Giroux, which includes the non-orientable case 
(see section \ref{egx}).
 An open book gives rise to a foliation $\mathcal O$ constructed as follows. 
 The pages endow 
 $B$ with a normal framing. So a tubular neighborhood $T$ of $B$ is 
trivialized:
 $T\cong B\times D^2$. Out of $T$ the leaves are the pages modified by spiraling
 around $T$; the boundary of $T$ is a union of compact leaves; and the 
interior of
 $T$ is foliated by a Reeb component, or a generalized Reeb component 
in the sense
 of Wood \cite{wood}. 
 For technical reasons in the homotopy argument of section \ref{homotopy}, the Reeb components of $\mathcal O$, instead of being usual Reeb components, will be 
 {\it thick Reeb components} in which 
 a neighborhood of the boundary is foliated by toric compact leaves.
 We call such a foliation an {\it open book foliation}. 
 
 The latter can be modified by inserting a so called 
{\it suspension foliation}. Precisely, let $\Si$ be a compact 
sub-surface of some leaf of $\mathcal O$ out of $T$ and 
 $\Si\times [-1,+1]$ be a foliated neighborhood of it (each  
$\Si\times \{t\}$ being contained in 
  a leaf of 
 $\mathcal O$). Let $\vp:\pi_1(\Si)\to {\rm Diff_c}(]-1,+1[)$
 be some 
 representation into the group of compactly supported diffeomorphisms;
$\vp$ is assumed to be trivial on the peripheral elements. It allows 
 us to construct a {\it suspension foliation} 
$\mathcal F_\vp$ on $\Si\times [-1,+1]$, whose leaves are transverse 
to the vertical segments $\{x\}\times [-1,+1]$ and whose holonomy
 is $\vp$. The modification consists of removing $\mathcal O$ from the
 interior of
 $\Si\times [-1,+1]$ and replacing it by $\mathcal F_\vp$. The new foliation,
 denoted 
 $\mathcal O_\vp$, is an {\it open book foliation modified by suspension.} 
 Theorem \ref{reg} can now be made more precise:
 
 \begin{thm}\label{modele} 
 Every co-orientable $\Ga_1^r$-structure, $r\geq 1$, is homotopic to an
 open book foliation modified by suspension.
 \end{thm}
 
 The proof of this theorem is given in sections \ref{tsub-c} - \ref{homotopy} when $r\geq 2$. 
In section \ref{BV}, %this theorem will be extended to $r\geq 1$ including the case
%$r=1+bv$.
we explain how to get the less regular case $1\leq r<2$. We have chosen to treat the case $r=1+bv$ (the holonomy local diffeomorphisms are $C^1$ and their first derivatives have
a bounded variation). Indeed, Mather observed in \cite{mather3} that 
$Diff^{1+bv}_c(\R)$ is not a perfect group and it is often believed that the perfectness of 
$Diff^r_c(\R)$ plays a role in the regularization theorem.

 %The case $r=1+bv$ will be treated in section 5.
 In section \ref{plane},  
 the homotopy class of the tangent plane field will be discussed. Finally
 the case of 
 non co-orientable $\Ga_1^r$-structure will be sketched in section
 \ref{comments} where the corresponding models, based on {\it twisted
 open book}, will be presented.

 We are very grateful to Vincent Colin,  \'Etienne Ghys and 
 Emmanuel  Giroux for  their comments, suggestions and explanations.

\vskip 1cm
\section{Tsuboi's construction }\label{tsub-c}

\medskip

%Every manifold $W$ carries the {\it trivial} (co-orientable) $\Ga_1$-structure, whose 
%leaves in $W\times \R$ are just the $W\times\{t\}$'s.
A $\Ga_1$-structure $\xi $ on $M$ is said to be {\it trivial}
 on a codimension 0 submanifold $W$  
when, for every $\vert t\vert$ small enough, $W\times\{t\}$ lies 
in a leaf of the associated foliation.

Every  closed 3-manifold $M$ has a so-called {\it Heegaard decomposition} 
$M=H_-\mathop{\cup}\limits_{\Si}H_+$, where $H_\pm$ is a possibly 
non-orientable handlebody (a ball with handles of index 1 attached) 
and $\Si$ is their common boundary. A {\it thick 
Heegaard decomposition} is a similar decomposition where the surface
 is thickened:
$$M=H'_-\mathop{\cup}\limits_{\Si\times\{-1\}} 
\Si\times[-1,+1] \mathop{\cup}\limits_{\Si\times\{+1\}}H'_+.
$$
The following statement  is due to T. Tsuboi in \cite{tsuboi} 
where it is left to the reader
as an exercise.

\begin{prop}\label{prereg}
 Given a $\Ga_1$-structure $\xi $ of class $C^r$, $r\geq 2$,
  on a closed 3-manifold $M$, there  
exists a thick Heegaard decomposition and a homotopy 
$\left(\xi_t\right)_{t\in[0,1]}$ from $\xi$ such that:

{ 1)}  $\xi_1$ is trivial on $H'_\pm$;

{ 2)} $\xi_1$ is regular on $\Si\times[-1,+1]$ and the induced foliation 
is a suspension.\\
\end{prop}

 \proof  With $\xi $ and its foliation $\mathcal F$ defined on an open 
  neighborhood of the zero section
 $M\times 0$ in $M\times \R$,  there comes a covering 
 of the zero section by boxes, open in $M\times\R$~, bi-foliated
  with respect to $\mathcal F$ and the fibers. We choose
 a $C^1$-triangulation  $Tr$ 
 of $M$ so fine that each simplex lies entirely in a box. With $Tr$
 comes a vector field
 $X$ defined as follows.
 
 First, on the standard $k$-simplex there is a smooth vector field $X_{\De^k}$,
  tangent to each face,
  which is the (descending) gradient  of a Morse function having one 
critical point 
 of index $k$ at the barycenter and one critical point of index $i$ at
 the barycenter of each $i$-face. When $\De^i\subset \De^k$
  is an $i$-face, $X_{\De^i}$ is the restriction of $X_{\De^k}$ to $\De^i$. 
Now, if
  $\si$ is a $k$-simplex of $Tr$, thought of 
  as a $C^1$-embedding $\si:\De^k\to M$, we 
  define $X_\si:=\si_*(X_{\De^k})$. The union of the $X_\si$'s is a 
$C^0$ vector field $X$
  which is uniquely integrable. After a reparametrization of each simplex
  we may assume that the stable manifold $W^s(b(\si))$ of the barycenter 
$b(\si)$
  is $C^1$.
  
  The $\Ga_1$-structure $\xi$ (co-oriented by the 
  $\R$ factor of $M\times\R$) is said to be in {\it Morse position}
 with respect to $Tr$ if:
  \bi
  \item[(i)] it has a smooth Morse type singularity of index $k$ at the 
barycenter of each $k$-simplex and it is regular elsewhere;
  \item[(ii)] $X$ is (negatively) transverse to $\xi $ out of the 
singularities.
  \ei
 
  \begin{lemme} \label{morsepos} Let $\mathcal F$
be the foliation associated to $\xi$. There exists 
  a smooth section $s$ such that 
   $s^*\mathcal F$ is in Morse position with respect to $Tr$.
  \end{lemme}
  
  Note that, as $s$ is homotopic to the zero section, the 
  $\Ga_1$-structure   $s^*\mathcal F$ on $M$ is homotopic to $\xi$.\\

  \proof Assume that  $s$ is already built near the $(k-1)$-skeleton. Let $\si$
  be a $k$-simplex. We explain how to extend $s$ on a neighborhood of 
  $\si$. After a fibered isotopy of $M\times\R$
  over the identity of $M$, we may assume that 
  $\mathcal F$  
  is trivial near $\{b(\si)\}\times\R$. Now, near
  $b(\si)$, we ask 
  $s$ to coincide with the graph of some local positive Morse 
  function $f_\si$ whose Hessian is negative definite on $T_{b(\si)}\si$ 
  and positive definite on $T_{b(\si)}W^s(b(\si))$. This function is now 
  fixed up to a positive 
  constant factor. We will extend $s$ as the graph of some function $h$
  in the $\mathcal F$-foliated chart over a neighborhood of $\si$. 
 This function is 
  already given on a neigborhood $N(\partial\si)$ of $\partial\si$ where 
  it is $C^r$, the regularity of $\xi$, and  satisfies
   $X.h<0$ except at the barycenter of each face.  
   
   On the one hand, choose an arbitrary extension $h_0$ of $h$ 
to a neighborhood of 
   $\si$ vanishing near $b(\si)$.
 On the other hand, choose a nonnegative function $g_\si$ such that:
   \bi
   \item[-] $g_\si=0$ near $\partial\si$;
   \item[-]  $g_\si=f_\si$ near $b(\si)$;
   \item[-] $X.g_\si<0$ when $g_\si>0$ except at $b(\si)$.
   \ei
   Then, if $c>0$ is a large enough constant, $h:=h_0+cg_\si$ has 
the required properties, except smoothness. Returning to $M\times\R$, 
the section $s$ we have
   built is $C^r$, smooth near the singularities, and $X $ is transverse 
to $s^*\xi$
   except at the singularities. Therefore, there exists a smooth 
$C^r$-approximation
   of $s$, relative to a neighborhood of the barycenters which meets all 
the required
   properties. \bull\\

 Thus, by a deformation of the zero section 
which induces a homotopy of $\xi$,
 we have  put $\xi$ in {\it Morse position} with respect to $Tr$. 
 In the same way, applying lemma \ref{morsepos} to the trivial
 $\Ga_1$-structure $\xi_0$, we also have a Morse function $f$ 
such that $X.f<0$ except at the barycenters.

Let $G_-$ (resp. $G_+$) denote the closure of the union of the unstable 
(resp. stable) manifolds of the singularities of $X$ of index 1 (resp. 2).
The following properties are clear: 
\bi
\item[(a)] In $M$, the subset $G_-$ (resp. $G_+$) is a $C^1$-complex of 
dimension 1.
\item[(b)] It admits arbitrarily small handlebody neighborhoods $H'_-$ 
(resp. $H'_+$) whose boundary is transverse to $X$. 
\item[(c)] Every orbit of $X$
outside $H'_\pm$ has one end point on $\partial H'_-$ and the other
 on $\partial 
H'_+$. This also holds true for any smooth $C^0$-approximation $\widetilde X$
of $X$ (in particular $\widetilde X$ is still negatively transverse to $\xi$).
\ei

Given a (co-orientable) $\Gamma_1$-structure $\xi$ on a space $G$, by an
 {\it upper (resp. lower) completion}
 of $\xi$ one means a foliation $\mathcal F$
of $G\times(-\epsilon, 1]$ (resp. $G\times[-1,\epsilon)$),
 for some positive $\epsilon$, which is transverse to every fiber
 $\{x\}\times(-\epsilon,1]$ (resp. $\{x\}\times[-1,\epsilon)$), 
whose germ along
$G\times \{0\}$ is $\xi$, and such that $G\times \{t\}$ is a leaf of
 $\mathcal F$ for every $t$ close enough to $+1$ (resp. $-1$).\\

\begin{lemme} \label{lower} 
Every co-orientable  $\Gamma_1^r$-structure on a simplicial complex $G$ of 
dimension 1, $r\ge 2$~, admits an upper (resp. lower)
completion of class $C^r$~.
\end{lemme}

 \nd {\sc Proof.}
One reduces immediately to the case where $G$ is a single edge.
In that case, using a partition of unity, one builds a line field which 
fulfills the claim.
 This line field is integrable.% since $G$ is of dimension 1.
 
${}$\bull\\

By (a), the $\Gamma_1$-structure
 $\xi$ admits an upper (resp. lower) completion over $G_+$
(resp. $G_-$), and thus also over an open neighborhood $N_+$ (resp. $N_-$)
of $G_+$ (resp. $G_-$). By (b), there is  
 a handlebody neighborhood $H'_\pm$ of $G_\pm$
contained in $N_\pm$ and whose boundary is transverse to $X$.
So we have a foliation $\mathcal F$ defined on a neighborhood 
of$$(M\times \{0\})\cup (H'_-\times[-1,0])\cup(H'_+\times[0,1])$$
which is transverse to $X$ on $M\smallsetminus (H'_-\cup H'_+)$ 
and tangent to  $H'_\pm\times \{t\}$ %is tangent to $\mathcal F$ 
for every $t$ close to $\pm 1$.

By (c), there is a diffeomorphism 
 $F: M\smallsetminus Int(H'_-\cup H'_+)\to \Sigma\times[-1,+1]$ for some
 closed surface $\Sigma$, which  maps orbit segments of $\widetilde X$ onto 
fibers.

For a small $\epsilon>0$, choose a function $\psi:{\R}\to[-1,+1]$
which is smooth, odd, and such that: 
\bi
\item $\psi(t)=0$ for $0\le t\le 1-3\epsilon$ and 
$\psi(1-2\epsilon)=\epsilon$; 
\item $\psi$ is affine on the interval $[1-2\epsilon,1-\epsilon]$; 
\item $\psi(1-\epsilon)=1-\epsilon$ and $\psi(t)=1$ for $t\ge 1$;
\item  $\psi'>0$ on the interval $]1-3\epsilon,1[$.
\ei

Let $s:M\to M\times\R$ be the graph of the function whose value is
 $\pm1$ on  $H'_\pm$ and  $\psi(t) $ at the point 
$F^{-1}(x,t)$ for $ (x,t)\in \Sigma\times[-1,+1]$. When  $\epsilon$ is small
 enough,
it is easily checked that, for every $x\in \Si$, the path 
$t\mapsto s\circ F^{-1}(x,t)$ is transverse to $\mathcal F$ except at its end 
points.
Then,  $\xi_1:=s^*\mathcal F$ is homotopic to $\xi$ and obviously fulfills the
 conditions 
required in  proposition \ref{prereg}. \bull\\
\medbreak

\vskip 1cm
\section{Giroux's construction }\label{egx}
\medskip
We use here theorem III.2.7 from Giroux's article \cite{giroux}, which 
states the following: {\it
\begin{itemize} 
\item[] Let $M$ be a closed 3-manifold (orientable or not). There exist a
 Morse function $f: M\to \R$ and a co-orientable surface $S$ which is 
$f$-essential in $M$. 
\end{itemize}
}
Giroux says that $S$ is {\it $f$-essential} when the restriction $f\vert S$ 
has exactly the same critical points as $f$ and the same local {\it extrema}. 
In the sequel, we call such a surface a {\it Giroux surface}.

Giroux explained to us 
\cite{giroux2} how this notion is related to  open book decompositions. 
In the above statement,
 the function $f$ can be easily  chosen {\it self-indexing} 
 (the value of a critical point is its Morse index in $M$). 
Thus, let $N$ be the level set 
 $f^{-1}(3/2)$. The smooth curve $B:=N\cap S$ will be the binding of the 
open book
 decomposition we are looking for. It can be proved that the following holds 
 for every  regular value $a, \ 0<a\leq 3/2$ :
 \begin{itemize}
 \item the level set $f^{-1}(a)$ is the union along their common boundaries 
of two surfaces, $N_1^a$ and $N_2^a$,
 each one being  
 diffeomorphic to the sub-level surface $S^a:= S\cap f^{-1}([0,a])$;
 \item the sub-level $M^a:= f^{-1}([0,a])$ is divided by $S^a$ into 
two parts $P_1^a$ and $P_2^a$ which are isomorphic handlebodies (with corners);
 \item $S^a$ is isotopic to $N_i^a$ through $P_i^a$, for $i=1,2$, by an 
isotopy fixing
 its boundary curve $S^a\cap f^{-1}(a)$.
 \end{itemize}
 This claim  is obvious when $a$ is small 
 and the property is preserved when crossing the critical level 1. 
 In this way the handlebody
 $H_-:= f^{-1}([0,3/2])$  is divided by $S^{3/2}$ into two diffeomorphic 
 parts $P_i^{3/2}$,
 $i=1,2$, and we have $N= N_1^{3/2}\cup N_2^{3/2}$. % where $N_i:= N\cap H_-^i$. 
 We take
 $S^{3/2}$, which is isotopic to $N_i^{3/2}$ in $P_i^{3/2}$, as a page. 
The figure is the same in 
 $H_+:=f^{-1}([3/2, 3])$. The open book decomposition is now clear. \\
 
 \begin{prop}\label{gir} 
 Let $K\subset M$ be a compact connected co-orientable surface 
 whose boundary is not empty. Then there exists an open book decomposition 
 whose
 some page contains $K$ in its interior.
 \end{prop}
 
 \proof (Giroux) According to the above discussion it is sufficient to
 find a Morse function $f$
and a Giroux  surface $S$ (with respect to $f$) containing $K$. 
Let $H_0$ be the quotient  of $K\times [-1,1]$ by shrinking to a point each
 interval
$\{x\}\times [-1,1] $ when $x\in \partial K\times [-1,1] $. 
After smoothing, it is a handlebody whose boundary is the double of $K$. 
On $H_0$ there exists a standard Morse function $f_0$ which 
is constant on $\partial H_0$, having one minimum, 
the other critical points being of 
 index 1. The surface $K\times\{0\}$ can be made $f_0$-essential.
 This function is then extended to a global Morse function $\tilde f_0$ on $M$.
 At this point we have to follow the proof of Theorem III.2.7 in \cite{giroux}.
 The function  $\tilde f_0$ is changed  on the complement of $H_0$, step by 
step when crossing its critical level, so that $K\times\{0\}$ extends as a 
Giroux surface in $M$.
 \bull\\
 
 Let now $\xi$ be a $\Ga_1$-structure meeting the conclusion of proposition
 \ref{prereg}, up to rescaling the interval to $[-\ep,+\ep]$. 
Let $\mathcal F_\vp$ be the 
suspension foliation induced by $\xi$
 on $\Si\times[-\ep,+\ep]$.
 Choose $x_0\in \Si\times\{0\}$; the segment $x_0\times[-\ep,+\ep]$
 is transverse to %the suspension foliation, 
 $\mathcal F_\vp$. Let $K$ be the surface obtained from $\Si\times 0$
 by removing a small open disk centered at $x_0$.  The foliation 
 $\mathcal F_\vp$ foliates $K\times[-\ep,+\ep]$ so that $K\times \{t\}$ 
 lies in a leaf, when $t$ is close to $\pm\ep$,
 and $\partial K\times[-\ep,+\ep]$ is foliated by parallel circles. We apply 
 proposition \ref{gir} to this $K$.

 \begin{cor} There exists an open book foliation $\mathcal O$ of $M$ inducing 
 the trivial foliation on $K\times[-\ep,+\ep]$
  (the leaves are $K\times\{t\},\ t\in [-\ep,+\ep] $).
 \end{cor}
 
 Therefore, we have an open book foliation modified by suspension by 
replacing the 
 above trivial foliation of
 $K\times[-\ep,+\ep]$ by $\widetilde{\mathcal F}_\vp$, the trace
 of $\mathcal F_\vp$ on $K\times[-\ep,+\ep]$. Let  $\mathcal O_\vp $
 be the resulting foliation of $M$ and $\xi_\vp$ be its regular
 $\Ga_1$-structure. For proving theorem \ref{modele} (when $r\geq 2$)
it is sufficient to prove that $\xi $ and $\xi_\vp$ are homotopic. 
This is done in the next
section. 

\vskip 1cm
\section{Homotopy of $\Ga_1$-structures}\label{homotopy}
\medskip
We are going to describe a homotopy from $\xi_\vp$ to $\xi $. Recall
the tube $T$ around the binding. For simplicity, we assume that 
 each component of $T$ is foliated by a standard
Reeb foliation; the same holds true if $T$ is foliated by Wood components (in the sense 
of \cite{wood}).  Let $T'$ be a slightly larger tube.   

\begin{lemme} \label{reeb}
There exists a homotopy, relative to $M\smallsetminus int(T')$,
 from $\xi_\vp$ to a new $\Ga_1$-structure $\xi_1$  on $M$ such that:
 
 \nd 1) $\xi_1$ is trivial on $T$;
 
 \nd 2) $\xi_1$ is regular on $T'\smallsetminus int(T)$ with compact  
toric leaves 
 near $\partial T$ and  spiraling half-cylinder  leaves with boundary in 
$\partial T'$ (as in an open book foliation).
 \end{lemme}
 
 \proof Recall from the introduction that we only use {\it thick} Reeb components. So there is a third concentric tube $T''$, $T\subset int(T'')\subset int(T')$, so that 
 $T''\smallsetminus int(T)$
  is foliated by toric leaves; and $int(T)$ is foliated by planes.

 On each component of $\partial T$ we have coodinates $(x,y)$
 coming from the
 framing of the binding, the $x$-axis being a parallel and the $y$-axis
 being a meridian.
 Let $\ga$ be a parallel in $\partial T$. The $\Ga_1$-structure 
which is induced by 
 $\xi_\vp$ on $\ga$ is singular but not trivial and the germ of 
foliation $\mathcal G$
 along the zero section 
 in the normal line bundle $A\cong \ga\times\R$ is shown on figure 1.
 
The annulus $A$ is endowed with coordinates $(x,z)\in\ga\times\R$. 
The orientation of 
the $z$-axis, which is also the orientation of the normal bundle to 
the foliation 
$\mathcal O_\vp$ along $\ga$, points to the interior of $T$. So the leaves of 
$\mathcal G$ are parallel circles in $\{z\leq 0\}$ and spiraling leaves 
in $\{z>0\}$. Take coordinates $(x,y,r)$ on $T$ %(also on $T''$) 
where $r$ is the distance
 to the binding; say that $r=1$ on 
$\partial T''$ and $r=1/2$ on $\partial T$. 
%For $\ep>0$, consider the map $g_\ep: T\to A$, $g_\ep (x,y,r)= (x, -\ep r^2+\ep) $.
Let $\la(r^2)$ be an even smooth function with $r=0$ as unique critical point,
vanishing at $r=1/2$ and 
$\la(1)<0<\la(0)$. Consider $g: T''\to A$, $g(x,y,r)= \bigl(x, \la(r^2)\bigr) $.

 \bigskip
\begin{center}
\begin{picture}(0,0)%
\epsfig{file=figure1.pstex}%
\end{picture}%
\setlength{\unitlength}{1579sp}%
\begingroup\makeatletter\ifx\SetFigFont\undefined%
\gdef\SetFigFont#1#2#3#4#5{%
  \reset@font\fontsize{#1}{#2pt}%
  \fontfamily{#3}\fontseries{#4}\fontshape{#5}%
  \selectfont}%
\fi\endgroup%
\begin{picture}(7212,6249)(3404,-6673)
\put(9601,-4111){\makebox(0,0)[lb]{\smash{{\SetFigFont{11}{13.2}{\familydefault}{\mddefault}{\updefault}{\color[rgb]{0,0,0}$z=0$}%
}}}}
\end{picture}%

\vskip .5cm
 Figure 1\\

\end{center}
\bigskip
It is easily seen
%that the Reeb foliation on $int(T'')$ is the pullback $g^*\mathcal G$; even, we have:
that $\xi_\vp\vert int(T'')\cong g^*\mathcal G$. 
%We let $\ep $ 
%go to 0, which is a homotopy.
%As the holonomy diffeomorphisms of $\mathcal G$ are $C^\infty$-tangent to 
%Identity 
 %along the compact leaf $\{z=0\}$, it makes sense to glue together 
%along $\partial T$
 %both  $\Ga_1$-structures: $\xi_\vp$ restricted to $M\smallsetminus int(T)$ and 
 %$g_0^*\mathcal G$ on $T$. 
 Let now $\bar\la(r^2)$ be a new even function coinciding with $\la(r^2)$ near
 $r=1$, having negative values everywhere and whose critical set is $\{r\in [0,1/2]\}$. Let 
  $\bar g: T''\to A$, $\bar g(x,y,r)= \bigl(x, \bar\la(r^2)\bigr) $. A barycentric combination  of $\la$ and $\bar\la$ yields a homotopy from $g$ to $\bar g$ 
  which is relative to a neighborhood of $\partial T''$.
  The $\Ga_1$-structure $\xi_1$ we are looking for is defined by 
  $\xi_1\vert int(T'')= \bar g^*(\mathcal G)$ and $\xi_1=\xi_\vp $ on a neighborhood of 
  $M\smallsetminus int(T'')$.\bull\\
 
% This $\Ga_1$-structure is not trivial on $T$.
% But if the zero section of its normal bundle is slightly pushed into the
 %negative $z$-direction over $T$ and if the pushing vanishes out of $T'$,
 %we get the $\Ga_1$-structure we are looking for. \bull\\
 
 Recall the domain $K\times [-\ep,+\ep]$ from the previous section.
  After the following lemma we are done with the homotopy problem.
  %(given by lemma \ref{reeb})
  
 \begin{lemme} There exits a homotopy from $\xi_1$  to 
 $\xi $ relative to $K\times [-\ep,+\ep]$.
 \end{lemme}
 
 \proof Let us denote $M':= M\smallsetminus int(K\times [-\ep,+\ep])$ 
 which is a manifold with boundary and corners. It is equivalent to prove that the 
 restrictions of $\xi_1$ and $\xi$ to $M'$ are homotopic relatively to $\partial M'$.
 Consider the standard closed 2-disk $D=D^2$ endowed with the $\Ga_1$-structure 
 $\xi_D$ which is shown on figure 2.
 
  \bigskip
\begin{center}
\begin{picture}(0,0)%
\includegraphics{figure2.pstex}%
\end{picture}%
\setlength{\unitlength}{1737sp}%
\begingroup\makeatletter\ifx\SetFigFont\undefined%
\gdef\SetFigFont#1#2#3#4#5{%
  \reset@font\fontsize{#1}{#2pt}%
  \fontfamily{#3}\fontseries{#4}\fontshape{#5}%
  \selectfont}%
\fi\endgroup%
\begin{picture}(7887,6470)(3077,-6956)
\put(6076,-3886){\makebox(0,0)[lb]{\smash{{\SetFigFont{11}{13.2}{\familydefault}{\mddefault}{\updefault}{\color[rgb]{0,0,0}$d$}%
}}}}
\put(9826,-3811){\makebox(0,0)[lb]{\smash{{\SetFigFont{11}{13.2}{\familydefault}{\mddefault}{\updefault}{\color[rgb]{0,0,0}$Im(i)$}%
}}}}
\end{picture}%
%pstex_t
\vskip .5cm
Figure 2
\end{center}

\bigskip
It is trivial on the small disk $d$ and regular on the annulus $D\smallsetminus
int(d)$. In the regular part,  the leaves are circles near $\partial d$ and the other leaves are spiraling, crossing $\partial D$ transversely. One checks that the restrition of $\xi_1$ to $M'$ 
has the form $f^*\xi_D$ from some map $f:M'\to D$. We take $f\vert T: T\to d$ to be
 the open book
trivialization of $T$ (recall the binding has  a canonical  framing); $f\vert \partial K\times [-\ep,+\ep]$ to be the projection $pr_2$
onto $[-\ep,+\ep]$ composed with an embedding $i:[-\ep,+\ep]\to \partial D$ and 
$f$ maps each
 leaf of the regular part of $\xi_1$ to a leaf of the regular part of $\xi_D$. As $K$
does not approach $T$, we can take  $f(\partial M') = i([-\ep,+\ep])$; actually, except
near $T$, $f$ is given by the fibration over $S^1=\partial D^2$ of the open book decomposition.

%We observe 
Once $\xi$ has Tsuboi's form (according to proposition \ref{prereg}),
 the restrition of $\xi $ to $M'$ has a similar form: $\xi=k^*\xi_D$
for some map $k:M'\to D^2$. Recall that $M'$ is the union of two handlebodies
and a solid cylinder $D^2\times [-\ep,+\ep]$. Take $k$ to be $i\circ pr_2$ 
on the cylinder
and $k$ to be constant on each handlebody. Observe that $f$ and $k $ 
coincide on $\partial M'$. As $D$ retracts by deformation onto the image of $i$, one deduces that $f$ and $k$ are homotopic relatively to $\partial M'.$
\bull\\

This finishes the proof of theorem \ref{modele} when $r\geq 2$.

\medskip
\vskip 1cm
\section{The case $C^{1+bv}$ }\label{BV}
\medskip
 A co-oriented $\Ga^r_1$-structure  $\xi$ on $M$ can be realized by a foliation
 $\mathcal F$ defined  on a neighborhood of the 0-section in $M\times\R$; 
it is 
made of 
bi-foliated charts which are $C^\infty$ in the direction of the
 leaves and $C^r$ in
the direction of the fibers. Consider such a box $\mathcal U$ over an open
disk $D$ centered at $x_0\in M$; its trace on the $x_0$-fiber is an interval
$I$.
 Each leaf of $\mathcal U$ reads $z=f(x,t)$, $x\in D$, for some $t\in I$.
Here $f$ is a function which is 
smooth in $x$ and $C^r$ in $t$,  with $f(x_0,t)=t$; the foliating property is 
equivalent to  
$\frac{\partial f}{\partial t}>0$.
When $r= 1+bv$, there is a positive measure $\mu(x,t)$ on $I$, 
%$I_x:=\mathcal U\cap\{x\}\times\R$, 
without atoms and depending smoothly on $x$, such that:
$$(*)\quad\quad
\frac{\partial }{\partial t}f(x,t)-\frac{\partial }
{\partial t}f(x,t_0)= \int_{t_0}^t\mu(x,t).$$

\begin{prop}\label{bv} Theorem \ref{modele} holds true for any class of
 regularity 
$r\geq 1$ including the class $r=1+bv$.
\end{prop}

\proof The only part of the proof which requires some care of  regularity
 is section
\ref{tsub-c}, especially the proof of lemma \ref{lower}. Indeed, we have to 
avoid integrating 
$C^0$ vector fields.
 For proving lemma \ref{lower} with weak regularity we use the lemmas below
 which we shall 
 prove in 
 the case $r=1+bv$ only.\\

%Given a family $f_a: [0,1]\to \R$  of $C^r$ functions, $r\geq 1$, depending  $C^r$ on 
%the  parameter $a$ in an interval $I$, it is said to be {\it foliating} when the 
%$\frac{\partial}{\partial a}f_a(x)>0$ for every $(a,x)\in I\times[0,1]$. This condition %amounts to saying that the graphs of the $f_a$'s are the leaves of a foliation 
%(transversely  $C^r$) of some domain in the plane.

\begin{lemme}\label{compl}
 Let $f:D\times I\to \R$ be a $C^r$-function as above. Assume 
$I=]-2\ep,+\ep[$ and $f(x,-2\ep)>-1$ for every $x\in D$. Then there exists a 
function $F:D\times[-1,0]\to\R $ of class  $C^r$ such that:
\bi
\item[1)] $F(x,t)= f(x,t) $ when $t\in [-\ep, 0]$,
\item[2)] $F(x,t)=t$ when $t$ is close to $-1$,
\item[3)] ${\ds \frac{\partial F}{\partial t}>0}$.
\ei
\end{lemme}

\proof  Let $\mu(x,t)$ be the positive measure whose support is  $[-\ep, 0]$ 
such that formula $(*)$ holds for every 
$(x,t) \in D\times[-\ep, 0]$ and  $t_0=0$.
There exists another positive measure $\nu(x,t)$, smooth in $x$ and 
 whose support is contained in 
$]-1,-\ep]$, such that 
$$(**)\quad\quad f(x,-\ep)= -1+\int_{-1}^{-\ep} \left(\int_{-1}^t \nu(x,\tau)\right)dt.$$
Then a solution is 
$$F(x,t)=  -1+\int_{-1}^{t} \left(\int_{-1}^s
\bigl(\mu(x,\tau)+ \nu(x,\tau)\bigr)\right)ds.$$\bull
 
\begin{lemme}\label{relative} 
Let $A_1$ and $A_2$ two disjoint closed sub-disks of $D$. Let $F_1$ and 
$F_2$ be two solutions of lemma \ref{compl}. Then there exists a third solution 
which equals $F_1$ when $x\in A_1$ and $F_2$ when $x\in A_2$.
\end{lemme}

\proof  Both solutions $F_1$ and $F_2$ differ by the choice of the measure
$\nu(x,t)$ in formula $(**)$, which is $\nu_i$ for $F_i$. 
 Choose a partition of unity
$1=\la_1(x)+\la_2(x)$ with $\la_i=1$  on $A_i$. 
Then $\nu(x,t)= \la_1(x)\nu_1(x,t)
+\la_2(x)\nu_2(x,t)$ yields the desired solution.
\bull\\

The proof of proposition \ref{bv} is now easy. As already said, 
it is sufficient to prove lemma \ref{lower}
%the claim 1 from proposition \ref{prereg}
 in class $C^r$, $r\geq 1$.  
It is an extension 
problem of a foliation given near  the 1-skeleton $Tr^{[1]}\times\{0\}$ to
 $Tr^{[1]}\times[-1,0]$. One covers  $Tr^{[1]}$ 
by finitely many $n$-disks $D_j$. 
The problem is solved in each $D_j\times[-1,0]$ by applying lemma \ref{compl}.
 By applying
 lemma \ref{relative}  one makes the different extensions match together.\bull

\medskip
\vskip 1cm
\section{Homotopy class of plane fields}\label{plane}
\medskip
It is possible to enhance theorem \ref{reg} by prescribing
the homotopy class of the underlying  co-oriented plane field 
(see proposition \ref{nu} below). 
The question of doing the same with respect to theorem \ref{modele} is
 more subtil 
(see proposition \ref{nu2}).\\

\begin{prop} \label{nu}
Given a co-oriented $\Ga_1$-structure $\xi$ on the closed 3-manifold $M$
and a homotopy class $[\nu]$ of co-oriented plane field  in the tangent 
space $\tau M$,
there exists a  regular $\Ga_1$-structure 
$\xi_{reg}$ homotopic to $\xi$ whose underlying foliation 
$\mathcal F_{reg}\cap M$
 has a  tangent co-oriented plane field in the class $[\nu]$.
\end{prop}

Before proving it we first recall some well-known facts 
on %homotopy classes of 
co-oriented plane fields (see \cite{geiges}). 
 Given a base plane field $\nu_0$, a suitable Thom-Pontryagin 
construction yields a natural bijection between  the set of 
 homotopy classes of plane fields on $M$ 
and $\Om^{\nu_0}_1(M)$, 
the group of (co)bordism classes of $\nu_0$-framed and oriented
%1- dimensional closed submanifolds
closed (maybe non-connected) curves  in $M$.  
A $\nu_0$-framing of the curve $\ga$ is an isomorphism of fiber bundles
$\ep: \nu(\ga, M)\to \nu_0\vert\ga$, whose source 
is the normal bundle to $\ga$ in $M$. 
We denote $\ga^\ep$ the curve endowed with this framing.
Moreover, given $\ga^\ep$, if $\ga'$ is homologous to $\ga$ in $M$ there exists 
a $\nu_0$-framing $\ep'$
such that $(\ga')^{\ep'}$ is cobordant to $\ga^\ep$.\\

\nd {\bf Proof of \ref{nu}.} We can start with an open book 
foliation $\mathcal O_\vp$
yielded by  theorem \ref{modele}. Let $\nu_0$ be its tangent plane field. 
%(the modification by suspension does not affect its homotopy class). 
Near the binding,
 the meridian loops (out of $T$) are transverse to $\mathcal O_\vp$ and 
homotopic to zero in $M$. As a consequence, each
1-homology class  may be represented as well by a (multi)-curve in a page
 or by a connected curve out of $T$ positively transverse to all pages.
 We do the second
choice for  $\ga^\ep$, the $\nu_0$-framed curve whose cobordism class
encodes $[\nu]$ with respect to $\nu_0$.

Hence we are allowed to {\it turbulize} $\mathcal O_\vp$ along $\ga$. 
In a small tube 
$T(\ga)$ about $\ga$, we put %a Reeb component or 
a Wood component.
 Outside, the 
leaves are spiraling around $\partial T(\ga)$. Let $\mathcal O_\vp^{turb}$
be the resulting foliation. Whatever the chosen type of Wood component is, the 
$\Ga_1$-structures of $\mathcal O_\vp^{turb}$ and
$\mathcal O_\vp$ are homotopic by arguing as in section \ref{homotopy}.
 But the framing $\ep$ tells us which sort of Wood 
component will be convenient for getting the desired class $[\nu]$ 
(see lemma 6.1 in \cite{wood}).
\bull\\

In the previous statement, we have lost the nice model we found in theorem
\ref{modele}. Actually, thanks to a lemma of Vincent Colin \cite{colin2},
 it is possible to recover our model, at least when $M$ is orientable 
(see below proposition \ref{nu2}).

\begin{lemme} {\rm (Colin)} Let $(B,p)$ be an open book decomposition of 
$M$ and 
$\ga$ be a simple connected curve in some page $P$. Assume $\ga$ is
 orientation preserving. Then 
there exist a positive stabilization $(B', p')$ of $(B,p)$
 and a curve $\ga'$ in 
 $B'$ which is isotopic to $\ga$ in $M$. When $\ga$ is a multi-curve, 
the same holds true after a sequence of stabilizations.
\end{lemme}

The {\it positive Hopf open book decomposition} 
of the 3-sphere is the one whose binding
is made of two unknots with linking number +1; 
a page is an annulus foliated by fibers of the Hopf fibration $S^3\to S^2$. 
A positive stabilization is a ``connected sum'' with this 
open book. The new page $P'$ is obtained from $P$ by {\it plumbing} 
 an annulus $A$
whose core bounds a disk in $M$ (see \cite{giroux-good} for more details 
and other references).\\

\proof If $\ga$ is connected, only one stabilization is needed.
 We are going to explain this case only. A tubular neighborhood of $\ga$ in 
$P$ is an annulus.

Choose a simple arc $\al$ in $P$ joining $\ga$ to some component $\beta$ 
 of $B$ without crossing $\ga$ again. Let $\tilde\ga$ be a simple 
arc from $\beta$
 to itself which follows $\al^{-1}*\ga*\al$.  The orientation assumption 
implies that
 the surgery of $\beta$ by $\tilde\ga$  in $P$ provides a curve with  
two connected 
  components, one of them being isotopic to $\ga$ in $P$. Let $P_\pi$  
be the page opposite to $P$ and $R:P\to P_\pi$ 
  the time $\pi$ of a flow transverse to the pages (and stationary on $B$). 
  The core curve $C$ of 
  the annulus $A$ that  we use for the plumbing is the
 union $\tilde\ga\cup R(\tilde\ga)$.
  And $A$ is $(+1)$-twisted  around $C$ (with respect to its unknot framing) 
as in the Hopf 
  open book. Let $H$ be the 1-handle which is the closure
 of $A\smallsetminus P$.
  Surgering $B$ by $H$ provides the new binding.
 By construction, one of its components is isotopic to $\ga$. \bull\\

\begin{prop}\label{nu2} Let $\mathcal O_\vp$ be an open book foliation
 modified by 
suspension, whose its underlying open book is denoted $(B,p)$. 
Let  $\nu_0$ be its tangent co-oriented plane field. Let $\ga^\ep$
 be a $\nu_0$-framed
curve in $M$ and $[\nu]$ be its associated class of plane field. 
Assume 
$\ga$ is orientation preserving. % (that is, $w_1([\ga])=0$). 
Then there exists an open book foliation $\mathcal O'_\vp$ with 
the following properties:

1) its tangent plane field is in the class $[\nu]$;
%\vfill\eject

2) the suspension modification is the same %(in the same support $K\times[-\ep,+\ep]$)
for  $\mathcal O'_\vp$ as for $\mathcal O_\vp$ and 
is supported in $K\times[-\ep,+\ep]$;

3) as $\Ga_1$-structures, $\mathcal O_\vp$ and $\mathcal O'_\vp$ are homotopic.
\end{prop}

\proof As said in the proof of \ref{nu}, up to framed cobordism,
 $\ga^\ep$ may be chosen as a simple (multi)-curve 
in one page $P$ of $(B,p)$.  Applying 
Colin's lemma provides a stabilization $(B',p')$
 such that, up to isotopy,  $\ga$
lies in the new binding. Observe that, if $K$ is in $P$, $K $ 
is still in the new page $P'$;
hence 2) holds for any open book foliation carried by $(B',p')$.
 Once $\ga^\ep$ is
 in the binding, for a suitable Wood component foliating a
 tube about $\ga^\ep$, item 1) is fulfilled. Finally item 3) 
follows from item 2) and the proofs in section \ref{homotopy}.
 \bull\\

%\vskip 1cm\vfill\eject
\section{Case of a $\Ga_1$-structure 
with a twisted normal bundle}\label{comments}
\medskip

 What happens when  the bundle $\nu$
 normal to
$\xi$ is twisted? It is known that a necessary condition to 
regularization is the existence of a fibered embedding
 $i: \nu\to \tau M$ into the tangent fiber bundle to $M$. 
Conversely, assuming that  this condition is  fulfilled,
 we are going to state a normal form theorem analogous
 to theorem \ref{modele}. Since no step of the  previous proof 
can be immediately adapted to  this situation,
 we believe that it deserves a sketch of proof.\\

\begin{rien} \label{71}{\rm In the first step (Tsuboi's construction), 
we do not have ``Morse position'' with respect to 
a triangulation, since index and co-index of a singularity cannot
 be distinguished. Instead of lemma \ref{morsepos},
 we have the following statement.
\bi
\item[] {\it After some homotopy, $\xi$ has  Morse singularities and admits a 
pseudo-gradient whose dynamics has no recurrence 
(that is, every orbit has a finite 
length).}
\ei
Here, by a pseudo-gradient, it is meant 
 a smooth
section $X$ of $Hom(\nu, \tau M)$, a {\it twisted vector field} indeed, 
such that 
$X\cdot\xi<0$ except at the singularities 
(this sign is well-defined whatever a local  orientation of $\nu$, 
or co-orientation of $\xi$, is chosen); 
 such a pseudo-gradient 
always exists by using an auxiliary Riemannian metric.}
\end{rien}

\nd {\bf Sketch of proof.} Generically $\xi$ has Morse singularities. 
Let $X_0$ be a first pseudo-gradient, which is required to be the usual negative 
gradient in Morse coordinates near each singularity. Finitely many mutually disjoint 
2-disks of $M$ are chosen in regular leaves
of $\xi$
%transverse to $X_0$ 
such that every orbit of $X_0$ crosses  one of them.
 Following  Wilson's idea \cite{wilson},   
$\xi $ and $X_0$ are changed in a neighborhood $D^2\times [-1,+1]$ of
 each disk into a {\it plug} such that every orbit of the modified
 pseudo-gradient $X$
is trapped by one of the plugs.
 The plug has the mirror symmetry with respect to
$D^2\times \{0\}$. In $D^2\times [0,1]$ we just modify $\xi$ by 
introducing a cancelling pair of singularities, center-saddle. \bull\\

Let $G$ be the closure of the one-dimensional invariant manifold 
of all saddles. 
It is a graph.  We claim: {\it $\nu\vert G$ is orientable}. 
Indeed, we orient each edge from its saddle end point to its center end point. 
This is an orientation of
$\nu\vert G$ over the complement of the vertices.  
It is easily checked that this orientation
extends over the vertices. Thus $X$ becomes a usual vector field near 
$G$ and we have an arbitrarily small tubular neighborhood $H$ 
of $G $ whose boundary is transverse to 
$X$, and $X$ enters $H$. Now, the negative completion of $\xi\vert H$ 
can be performed  as in lemma \ref{lower}.

The complement  $\hat M$ of $int\,H$ in $M$ is fibered over a surface $\Si$, 
the fibers being intervals ($\cong [-1,1]$)
tangent to $X$. By taking a section we think of $\Si$ 
as a surface in $M\smallsetminus H$.  Since $\xi$ is not co-orientable,
 {\it $\Si$ is one-sided and $G$ is 
connected.} Arguing as in section \ref{tsub-c}, {\it after some homotopy,
 $\xi$ becomes trivial on $H$ and transverse to $X$ on $\hat M$,
 hence a suspension foliation corresponding to a
 representation $\vp:\pi_1(\Si)\to Diff_c(]-1,1[)$. }\\

\begin{rien} \label{72}{\rm 
In the second step (Giroux's construction), we have to leave the open books
 and we need a {\it twisted open book}. It is made of the following:
\bi
\item a binding $B$ which is a 1-dimensional closed co-orientable submanifold
 in $M$;
\item a {\it Seifert fibration} $p: M\smallsetminus B\to [-1,+1]$ which has  two
 one-sided  exceptional surface
fibers $p^{-1}(\pm 1)$ and which is a proper smooth 
submersion over the open interval;
\item when $t$ goes to $\pm 1$, $p^{-1}(t)$
 tends to a 2-fold covering of $p^{-1}(\pm 1)$;
\item near $B$ the foliation looks like an open book.
\ei
%}\end{rien}

The exceptional fibers are compactified by $B$ as smooth surfaces with
 boundary. But,
for $t\in\,  ]-1,+1[$, $p^{-1}(t)$ is compactified by $B$ as a closed 
surface showing
(in general)
an angle along  $B$. Notice that, since $B$ is co-orientable,
 a twisted open book 
gives rise to a smooth foliation where each component of the binding is 
replaced by a Reeb 
component, the pages spiraling around it.
}\end{rien}

Such an open book is generated by a {\it one-sided Giroux surface}, 
which is the union of the compactified exceptional fibers. Abstractly,
 a one-sided Giroux surface in 
$M $ with respect to a 
Morse function $f:M\to \R$ is a one-sided surface $S$ such 
that $f\vert S$ has the same critical points and the same extrema as $f$  
and fulfills the extra condition: for every  regular value 
$t\in \R$, $f^{-1}(t)\cap S$ is a two-sided curve 
in the level set  $f^{-1}(t)$. Starting with $(S,f)$ where $f$ is a 
self-indexing Morse function, a twisted open book is easily constructed. 
Its binding is the co-orientable curve $f^{-1}(3/2)\cap S$. In general 
such a one-sided Giroux surface (or twisted open book) does not exist on $M$;
 the obstruction lies in the existence of a twisted line subbundle of
 $\tau M$. 
Nevertheless, with a suitable assumption, we have an 
analogue of proposition \ref{gir}:
\bi
\item[]{\it Let $i:\nu\to\tau M $ be an embedding of a twisted line bundle.
Let $K\subset M$ be a compact connected one-sided surface
 whose boundary 
is not empty. Assume the following: $\nu\vert K$ is twisted,
  $\nu\vert\partial K$ is trivial
 and the normal bundle $\nu(K,M)$ is homotopic 
to $i(\nu)\vert K$. Then there exists a twisted open book with one 
exceptional page
containing $K$ in its interior.}\\
\ei
\nd {\bf Sketch of proof.} We give the proof only in the setting of \ref{71} by taking 
$K$  to be the closure of $\Si\smallsetminus d_0$ where  $d_0$ is a small closed 
2-disk  in $\Sigma$. Recall the fibration
$\rho:\hat M\to \Si$. Let $H'$ be the handlebody made of  $H$ to which is glued
the 1-handle $\rho^{-1}(d_0)$.  Let $M'$ be the complement of $int\,H'$. 
Consider   a minimal system of mutually disjoint 
 compression disks $d_1,...,d_g$
 of $H$, so that cutting $H$ along them yields a ball; $g$  is the genus of $H$. 
 Thus $d_0, d_1,...,d_g$ is a minimal system of mutually disjoint compression disks of $H'$. 
 
 We claim: {\it $g+1$ is even.} Indeed, by assumption
  there exists a non vanishing section of the bundle
 $Hom(\nu,\tau M)$. Thus, the number
of zeroes of the pseudo-gradient $X$ is even and 
  the Euler characteristic $\chi(G)$ 
  of the graph $G$ is even. As $G$ is connected, the genus $g$ is odd, which proves
  the claim.
  
Now, one follows Giroux's algorithm for completing $K$ to a closed 
Giroux's surface. On the surface $\partial M'$ the union 
of the attaching curves $\partial d_0, \partial d_1,...,
 \partial d_g$ is not separating. Then, after some isotopy, for each $i=0,..., g$, 
 $\partial d_i$  crosses
 $\partial K$ in exactly two points $a_i,b_i$ linked by an arc $\al_i$ (resp. $\al'_i)$ 
 in  $\partial d_i$ (resp. $\partial K$), so that $\al_i\cup\al'_i$ bounds a disk 
 in $\partial M'$. Moreover, one can arrange that all the arcs $\al_0, ..., \al_g$ 
 are parallel. Also we link $a_i, b_i$ by a simple arc in $d_i$.
 Now, each compression disk defines,
 simultaneously,
 a 1-handle which  is glued 
 to $K$ and a 2-handle which is glued to $M'$, yielding a proper surface $K_1$ in some
 3-submanifold $M''$ of $M$, whose complement is a ball. The boundary of $K_1$
 is made of $g+2$  parallel curves in the sphere $\partial M''$. As this number is odd,
 Giroux described  a process of adding cancelling pairs of 1- and 2- handles whose effect is to change $K_1$ into $K_2\subset M''$ such that $\partial K_2$ is made of one curve only (\cite{giroux}, p. 676-677).
 Hence, $K_2$ can be closed into a Giroux's surface. \bull\\

 %We are ready to state the theorem.
 \begin{thm}
 Let $\xi$ be a non co-orientable $\Ga_1$-structure on $M^3$ 
whose normal bundle $\nu$ embeds into $\tau M$. Then 
 $\xi$ is homotopic to a twisted open book foliation modified by suspension.
 \end{thm}
 
 \nd {\bf Sketch of proof.} Let $K=\Si\smallsetminus int~D^2$ be the surface 
 with a hole, where $\Si$ was built in  the first step \ref{71}; 
it meets the required assumptions for building a twisted open book.
 
The twisted open book built in the second step gives rise to a foliation $\mathcal O$. 
Indeed, as the binding $B$ is co-orientable, it is 
 allowed to spiral the pages around a tubular neighborhood of $B$. 
The tube itself is foliated by thick Reeb components. As in the co-orientable case,
 we can modify  the open book foliation in a neighborhood of $K$ using
 the representation $\vp$, yielding the foliation $\mathcal O_\vp$
  and its associated regular
 $\Ga_1$-structure $\xi_\vp$. 
 %Arguing as in section \ref{homotopy}, 
%this foliation is homotopic to $\xi$.
We have to prove that $\xi$ and  $\xi_\vp$ are homotopic. We may suppose that 
$\xi $ is in Tsuboi form (\ref{71}). % Let $\xi_{Id}$ be the regular $\Ga_1$-structure %associated to the trivial representation.

We observe that the total space of $\nu$ has a foliation 
$\mathcal F_0$ (unique up to isomorphism) transverse to the fibers, having the zero section as a leaf and whose all non trivial holonomy elements
have order 2. It defines the trivial $\Ga_1$-structure $\xi_0$
in the twisted sense. Using 
notation of \ref{72}, one can prove that $\xi\vert H'$ and $\xi_\vp\vert H'$ are both homotopic to $\xi_0\vert H'$. Moreover, both homotopies coincide 
on the  boundary
$\partial H'$ (on $H'$, it is sufficient to think of the case when 
$\vp$ is the trivial representation. 
Thus  $\xi\vert H'$ and $\xi_\vp\vert H'$ are 
homotopic relative to the boundary. \bull\\

\begin{rien} Plane field homotopy class. By turbulizing $\mathcal O_\vp$, it it possible 
to have the normal field in any homotopy class of embeddings $\nu\to \tau M$.
\end{rien}
 
 Indeed, a curve in $M $ is homotopic to a curve transverse to $\mathcal O_\vp$ if and 
 only if it does not twist $\nu$. But these homology classes are exactly those which
 appear as a first difference homology class when comparing two embeddings $j_1,j_2
 :\nu\to\tau M$, since a closed curve which twits $\nu$ is not a cycle in 
 $H_1(M,\Z_{or(\nu)})\cong H^2(M, \Z_{or(\nu*\otimes\tau M)})$.

\end{document}